\documentclass[12pt]{article}
\usepackage{amsmath,amssymb,amsthm}
\usepackage{mathrsfs}  
\numberwithin{equation}{section}

\newtheorem{thm}{Theorem}[section]
\newtheorem{lem}{Lemma}[section]

\newtheorem{pro}{Proposition}[section]
\newtheorem{cor}{Corollary}[section]
\newcommand{\n}{\nonumber}

\newcommand{\w}{\omega}

\renewcommand{\l}{\lambda}

\renewcommand{\a}{\alpha}

\newcommand{\s}{\sigma}
\renewcommand{\o}{\omega}

\newcommand{\bb}{\begin{equation}}
\newcommand{\ee}{\end{equation}}
\newcommand{\bq}{\begin{eqnarray}}
\newcommand{\eq}{\end{eqnarray}}
\newcommand{\bqn}{\begin{eqnarray*}}
\newcommand{\eqn}{\end{eqnarray*}}

\newcommand{\R}{\mathbb{R}}

\newcommand{\N}{\mathbb{N}}

\newcommand{\bu}{v}

\newcommand{\bC}{\boldsymbol C}

\newcommand{\mes}{\operatorname{\rm meas}}

\newcommand{\const}{\operatorname*{const}}

\newcommand{\be}{\begin{equation}}
\newcommand{\bea}{\begin{eqnarray}}
\newcommand{\eea}{\end{eqnarray}}
\newcommand{\bean}{\begin{eqnarray*}}
\newcommand{\eean}{\end{eqnarray*}}

\newcommand{\var}{\varepsilon}
\renewcommand{\o}{\omega}

\newcommand{\intl}{\int\limits}

\newcommand{\Beweisende}{\rule{0.2cm}{0.2cm}}

\newcommand{\intmw}{{\int\hspace{-830000sp}-\!\!}}

\theoremstyle{definition}

\newtheorem{defin}[thm]{Definition}
\newtheorem{rem}[thm]{Remark}

\begin{document}
\title{On Liouville type theorems for the steady Navier-Stokes equations in $\Bbb R^3$}
\author{Dongho Chae$^*$  and J\"{o}rg Wolf $^\dagger$\\
\ \\
 $*$Department of Mathematics\\
Chung-Ang University\\
 Seoul 156-756, Republic of Korea\\
 e-mail: dchae@cau.ac.kr\\
and \\
$\dagger$Department of Mathematics\\
Humboldt University Berlin\\
Unter den Linden 6, 10099 Berlin, Germany\\
e-mail: jwolf@math.hu-berlin.de}
\date{}
\maketitle

\begin{abstract}
In this paper we prove three different Liouville type theorems for the steady Navier-Stokes equations in $\Bbb R^3$.
In the first theorem we improve logarithmically the well-known $L^{\frac92} (\Bbb R^3)$ result. In the second theorem we present a sufficient 
condition for the trivially of the solution($v=0$) in terms of the head pressure, $Q=\frac12 |v|^2 +p$. The imposed integrability condition here has  the same scaling property as the Dirichlet integral.
In the last theorem we present Fubini type condition, which guarantee $v=0$.
  \noindent{\bf keywords:} steady Navier-Stokes equations, Liouville type theorems
\end{abstract}
\section{Introduction}
\setcounter{equation}{0}
We consider the following stationary Navier-Stokes equations equations (NS) on $\mathbb{R}^3$.
$$
(NS)\left\{\aligned
(v\cdot\nabla)v&=-\nabla p + \Delta v,\\
\text{div}\, v &=0,
\endaligned
\right.
$$
where $v(x)=(v_1(x),v_2(x),v_3(x))$ and $p=p(x)$ for all $x\in\mathbb{R}^3$. The system is equipped with the boundary condition:
\begin{equation}
|v(x)|\rightarrow 0\quad \text{as} \quad|x|\rightarrow +\infty.\label{11}
\end{equation}
In addition to (\ref{11}) one usually also assume following finiteness of the  Dirichlet integral.
\begin{equation}
\int_{\mathbb{R}^3}|\nabla v|^2 dx < +\infty . \label{12}
\end{equation}
    A long standing open question is if any weak solution of (NS) satisfying the conditions \eqref{11} and \eqref{12}  is  trivial (namely, $v=0$ on $\mathbb{R}^3$).   We refer the book by Galdi(\cite{Galdi}) for the details on the motivations and historical backgrounds on the problem and the related results.
As a partial progress to the problem we mention that the condition $v\in L^{\frac{9}{2}}(\mathbb{R}^3)$ implies that
$v=0$ (see Theorem X.9.5, pp.729 \cite{Galdi}). As shown in \cite{cha1},  a different condition $\Delta v \in L^{\frac{6}{5}}(\mathbb{R}^3)$ also 
 imply $v=0$.  Another interesting progress, which shows that a solution $v\in BMO^{-1} (\Bbb R^3)$ to (NS), satisfying (\ref{12}) is trivial is obtained very recently by Seregin in \cite{ser}. For the case of plane flows the problem is solved by Gilbarg and Weinberger in \cite{Gil}, while the special case of the axially symmetric 3D flows without swirl is studied recently by Korobkov, M. Pileckas and R. Russo in \cite{Kor}(see also \cite{Koch}). In this paper we present three theorems, which present sufficient conditions to guarantee the triviality of the solution to (NS).

In the first theorem below we improve the above mentioned $L^{\frac92}$-result logarithmically.
\begin{thm}
\label{thm2.1}
Let $ \bu \in L^1_{ loc}(\R^{3})$ be a distributional solution to (NS) such that 
\begin{equation}
\intl_{ \R^{3}}  | \bu |^{\frac92}\left\{\log \Big(2+ \frac{1}{|\bu|}\Big)\right\}^{ -1} dx <+\infty.
\label{2.1aa}
\end{equation}   
Then $ \bu \equiv { 0}$.
\end{thm}

For discussion of the next theorem we introduce the head pressure,
$$Q=\frac12 |v|^2 +p,$$
which has an  important role in the study of the stationary Euler equations via the Bernoulli theorem. It is known(see e.g. Theorem X.5.1, pp. 688 \cite{Galdi}) that under the condition (\ref{11})-(\ref{12}) we have $p(x)\to p_0$ as $|x|\to +\infty$, where $p_0$ 
is a constant, which implies that  
\bb\label{13}
Q(x)\to 0 \qquad \mbox{as}\quad  |x|\to +\infty
\ee
after  re-defining $Q-p_0$ as the new head pressure.  Our second theorem below assumes integrability of $Q$ to conclude the triviality of $v$.
\begin{thm}
Let $(v,p)$ be a smooth solution to (NS) satisfying  (\ref{13}). Let us set $M:=\sup_{x\in \Bbb R^3} |Q(x)|$. 
 Then,  we have the  following inequality.
\bb\label{th11}
\int_{\Bbb R^3} \frac{|\nabla Q|^2}{|Q|} \left(\log \frac{eM}{|Q|} \right)^{-\a-1} dx\leq \frac{1}{\a} \int_{\Bbb \Bbb R^3} |\o|^2 dx \quad \forall \a >0.
\ee
Moreover, suppose there holds  the boundary conditions (\ref{11}), (\ref{13}) and 
\bb\label{th12}
\int_{\Bbb R^3} \frac{|\nabla Q|^2}{|Q|} \left(\log \frac{eM}{|Q|} \right)^{-\a-1} dx=o\left(\frac{1}{\a}\right)\quad \mbox{as}\quad \a \to 0,
\ee
then $v=0$  on $\Bbb R^3$.
\end{thm}
\noindent{\bf Remark 1.1.  } Since $|\nabla \sqrt{|Q|}|^2 = \frac14\frac{|\nabla Q|^2}{|Q|}$, and $\sqrt{|Q|}$ has the same scaling as the velocity
the integral $\int_{\Bbb R^3} \frac{|\nabla Q|^2}{|Q|}dx $ has the same scaling property as the Dirichlet integral in (\ref{12}).\\

Our third result concerns on the Fubini type condition for suitable function $\Phi(x,y)$ for $(x,y)\in \Bbb R^3\times \Bbb R^3$ to guarantee the triviality of 
the solution to (NS).

\begin{thm}
Let $v$ be a smooth solution to (NS) on $\mathbb{R}^3$ satisfying \eqref{11} and set $\o={\rm curl }\, v$.
Suppose there exists $q \in [\frac{3}{2},3)$ such that $x \in L^q (\mathbb{R}^3)$. We set
\begin{equation}
\Phi(x,y):= \frac{1}{4\pi} \frac{\w(x)\cdot(x-y)\times(v(y)\times \w(y))}{|x-y|^3}\label{1.5}
\end{equation}
for all $(x,y)\in \mathbb{R}^3 \times \mathbb{R}^3$ with $x\neq y$.
Then, it holds
\begin{equation}\label{2.19}
\int_{\mathbb{R}^3} |\Phi(x,y)|dy+\int_{\mathbb{R}^3} |\Phi(x,y)|dx <\infty \quad\forall(x,y)\in \mathbb{R}^3\times\mathbb{R}^3.
\end{equation}
Furthermore, if there holds 
\bb\label{th13}
\int_{\Bbb R^3}\int_{\Bbb R^3}  \Phi (x,y)dxdy=\int_{\Bbb R^3} \int_{\Bbb R^3} \Phi (x,y)dydx,
\ee
then,  $v=0$ on $\Bbb R^3$.
\end{thm}

\noindent{\bf Remark 1.2. } One can show that if $\w\in L^{\frac{9}{5}}(\mathbb{R}^3)$ is satisfied together with \eqref{11}, then
\eqref{th13} holds, and therefore $v$ is trivial. Although this result follows immediately by applying the $L^{\frac92}$-result together with Sobolev inequality
and the Calderon-Zygmund inequality, $\|v\|_{L^{\frac92}} \leq C \|\nabla v\|_{L^{\frac95}} \leq C \|\o\|_{L^{\frac95}}$. The above theorem provides us with different 
proof of this.   In order to check this  result we first recall the estimate of the Riesz potential on $\mathbb{R}^3$(\cite{Stein}),
\begin{equation}
\| I_{\alpha}(f)\|_{L^q}\leq C \| f \|_{L^p},\quad \frac{1}{q}=\frac{1}{p}-\frac{\a}{3},\quad
1\leq p < q < +\infty,\label{1.8}
\end{equation}
where
$$I_{\alpha}(f) := C \int_{\mathbb{R}^3}\frac{f(y)}{|x-y|^{3-\alpha}} dy,\quad 0<\alpha<3$$
for a positive constant $C=C(\alpha)$. Applying \eqref{1.8} with $\alpha=1$, we obtain by the H\"{o}lder inequality,
\begin{align*}
\int_{\mathbb{R}^3}\int_{\mathbb{R}^3} &|\Phi(x,y)| dydx
 \leq \int_{\mathbb{R}^3}\int_{\mathbb{R}^3} \frac{|\w(x)||\w(y)||v(y)|}{|x-y|^2} dydx\\
 &\leq \left(\int_{\mathbb{R}^3}|\w(x)|^{\frac{9}{5}}dx\right)^{\frac{5}{9}}
 \left\{ \int_{\mathbb{R}^3}\left(\int_{\mathbb{R}^3} \frac{|\w(y)||v(y)|}{|x-y|^2} dy\right)^{\frac{9}{4}}dx\right\}^{\frac{4}{9}}\\
 &\leq C\| \w\|_{L^{\frac{9}{5}}}\left(\int_{\mathbb{R}^3}|\w|^{\frac{9}{7}}|v|^{\frac{9}{7}}dx\right)^{\frac{7}{9}}\\
 &\leq C\| \w\|_{L^{\frac{9}{5}}}\left(\int_{\mathbb{R}^3} |\w|^{\frac{9}{5}}dx \right)^{\frac{5}{9}}
  \left(\int_{\mathbb{R}^3} |v|^{\frac{9}{2}} dx\right)^{\frac{2}{9}}\\
 &\leq C\|\w\|^2_{L^{\frac{9}{5}}}\|\nabla v\|_{L^{\frac{9}{5}}}\leq C\|\w\|^3_{L^{\frac{9}{5}}}<+\infty,
\end{align*}
Thus, by the Fubini-Tonelli theorem, \eqref{th13}  holds.

\section{ Proof of the main theorems}
\setcounter{equation}{0}

Below we use the notation $A \lesssim B$ if there exists an absolute constant $\kappa $ such that $A\leq \kappa B$.

\subsection {Proof of Theorem 1.1}

\begin{defin}
\label{def3.1}
Let $ \phi \in C^2(\R)$ be an N-function, i.\,e.  $ \phi $ is an even function such that 
$ \lim_{\tau  \to 0} \phi'(\tau) =0 $, and $ \lim_{\tau  \to \infty} \phi' (\tau)= +\infty$. We say $ \phi $ belongs 
to the class $ N(p_0, p_1)$ ($ 1< p_0 \le p_1 < +\infty$) if for all $ \tau  \ge 0$
\begin{equation}
(p_0-1)\phi'(\tau ) \le \tau \phi ''(\tau ) \le (p_1-1) \phi' (\tau ).
\label{3.1}
\end{equation}

\end{defin}

\begin{rem}
It is well known that $ \phi \in N(p_0, p_1)$ implies for all $ \tau \ge 0$
\begin{equation}
\phi(\tau ) \le \tau \phi '(\tau ) \le p_1\phi(\tau ).
\label{3.2}
\end{equation}
\end{rem}

\hspace{0.5cm}
We now define for   $ q>1$ 
\[
\phi _q(\tau ) =  \intl^{\tau }_{0} \frac{\xi ^{ q-1}}{ \log \frac{1+ 2\xi  }{\xi  }}  d\xi,\quad \tau \ge 0.  
\]
We easily calculate,
\begin{align*}
\phi_q '(\tau ) &= \frac{\tau ^{ q-1}}{ \log \frac{1+ 2\tau  }{\tau  }},\quad  
\\
\phi_q ''(\tau ) &= (q-1) \frac{\tau ^{ q-2}}{ \log \frac{1+ 2\tau }{\tau }} + \frac{\tau ^{ q-2}}{ \log^2 \frac{1+ 2\tau  }{\tau  }} \frac{1}{1+2\tau } 
\\
&= \frac{\phi_q '(\tau )}{\tau }\bigg(
(q-1) + \frac{1}{(1+2\tau ) \log \frac{1+2\tau }{\tau }}\bigg). 
\end{align*}
 Observing that $ \frac{1}{(1+2\tau ) \log \frac{1+2\tau }{\tau }} \le \frac{1}{\log 2}$, we get  for all $ \tau \ge 0$
\begin{equation}
(q-1) \phi_q '(\tau ) \le \tau \phi_q ''(\tau ) \le  (q+ (\log 2)^{ -1} -1)  \phi_q '(\tau ).
\label{3.3}
\end{equation}
This shows that $ \phi \in N(q, q+ (\log 2)^{ -1})$, and according to \eqref{3.2} it holds  
\begin{equation}
\phi_q(\tau ) \sim \tau \phi_q '(\tau ) = \frac{\tau ^{ q}}{ \log \frac{2+ \tau  }{\tau  }}. 
\label{3.4}
\end{equation}
Thus, \eqref{2.1} is equivalent to  
\begin{equation}
\intl_{ \R^3} \phi_{ \frac92} (| \bu |) dx < +\infty.   
\label{3.5}
\end{equation}

\begin{lem}
\label{lem3.3}
For any constant $ a > \frac{1}{2}$ we have 
\begin{equation}
\log a\frac{1+2\tau }{\tau } \sim \log \frac{1+2\tau }{\tau }.
\label{3.6}
\end{equation}

\end{lem}

\noindent{\bf Proof } In case $ a \ge 1$ we immediately get $\log a\frac{1+2\tau }{\tau } \ge  \log \frac{1+2\tau }{\tau } $. 
For the reverse we get for all $ 0 < \tau \le 1$, 
\[
\log a\frac{1+2\tau }{\tau } \le \Big(\frac{\log a}{\log 3} +1\Big) \log \frac{1+2\tau }{\tau },
\]
and for all $ \tau >1$
\[
\log a\frac{1+2\tau }{\tau } \le \log a + \log 3  \le \frac{\log a + \log 3 }{\log 2} \log \frac{1+2\tau }{\tau },
\]
which proves the claim. 

\hspace{0.5cm}
In case $ a <1$ we see that $\log a\frac{1+2\tau }{\tau } \le  \log \frac{1+2\tau }{\tau } $. 
On the other hand, we may choose $ \tau _0 >0$, such that
\[
\log \frac{1+2\tau_0 }{\tau_0 }= \frac{1}{2} \Big(1+ \frac{\log 2}{\log a^{ -1}}\Big)\log a^{ -1}.
\]
Then for $ \tau  \le  \tau _0$ we obtain 
\begin{align*}
\log a\frac{1+2\tau }{\tau } &= - \log a^{ -1} + \log \frac{1+2\tau }{\tau } 
= - 2 \Big(1+ \frac{\log 2}{\log a^{ -1}}\Big)^{ -1} \log \frac{2+\tau_0 }{\tau_0 } + \log \frac{1+2\tau }{\tau }. 
\\
&\ge \bigg[1- 2 \Big(1+ \frac{\log 2}{\log a^{ -1}}\Big)^{ -1}\bigg] \log \frac{1+2\tau }{\tau }
\\
& = \frac{\log 2- \log a^{ -1}}{\log 2 + \log a^{ -1}}\log \frac{1+2\tau }{\tau }. 
\end{align*}
For $ \tau > \tau _0$ we easily see that 
\begin{align*}
\log a\frac{1+2\tau }{\tau } & \ge \log 2 -   \log a^{ -1} = \frac{\log 2 -   \log a^{ -1} }{\log \frac{1+ 2\tau _0}{\tau _0}} 
\log \frac{1+2\tau_0 }{\tau_0 }
\\
&\ge \frac{\log 2 -   \log a^{ -1} }{\log \frac{1+ 2\tau _0}{\tau _0}} 
\log \frac{1+2\tau}{\tau}
\\
&=2 \frac{\log 2 -   \log a^{ -1} }{\log 2+ \log a^{ -1}} 
\log \frac{1+2\tau}{\tau}. 
\end{align*}
Whence, the claim.  \hfill \Beweisende 

\begin{lem}
\label{lem3.4} For all $ k\in \N$
\begin{equation}
\log \frac{1+2 \tau }{\tau } \sim \log \frac{1+ 2\tau ^k}{\tau ^k}, 
\label{3.7}
\end{equation}
where the hidden constants depend on $ q$ and $ k$ only. 

\end{lem}
\noindent{\bf Proof } In fact having  
$ 1+ 2\tau^k \le (1+2\tau )^k \le 2^{ k-1} (1+ 2^k\tau ^k) \le 2^{ 2k-2}(1+ 2\tau ^k)$ along with Lemma\,\ref{lem3.3}, we obtain 
\begin{align*}
\log \frac{1+ 2\tau^k }{\tau^k } \le k \log \frac{1+ 2\tau}{\tau} \le  \log 2^{ 2k -2}\frac{1+ 2\tau ^k}{\tau^k }
\lesssim  \log \frac{1+ 2\tau^k }{\tau^k }. 
\end{align*}
This proves the claim.  \hfill \Beweisende 

\begin{lem}
\label{lem2.3}
Let $ f \in L^1(\R^{3})$. Then for every $ \var >0$, there exists $ R > \var ^{ -1}$, such that 
\begin{equation}
\intl_{B_{ R}  \setminus B_{ R/2}} | f| dx \le  \frac{\var }{\log R}.  
\label{3.9a}
\end{equation}

\end{lem}

\noindent{\bf Proof  }  Assume the assertion of the lemma is not true. Then there exists $ \var >0$  such that  for all $ R \ge \var ^{ -1}$
\eqref{3.9a} does not hold. This implies for all  
 $ k \ge N$ with $ 2^k \ge \var ^{ -1}$
\[
\intl_{B_{ 2^k}  \setminus B_{ 2^{ k-1}}} | f| dx\ge  \frac{\var }{k\log 2}. 
\]
However the sum of  right-hand side from $ k=N$ to $ \infty$ is  infinite  which clearly contradicts to $ f\in L^1(\R^{3})$. 
Thus, the assumption is not true and therefore the assertion of the lemma holds.  \hfill \Beweisende      \\


 In view of \eqref{2.1aa} we easily see that $ \bu \in L ^{ \frac92}_{ \rm loc}(\R^{3})$. By using a standard mollifying argument we  
verify that $\bu \in W^{1,\, 2}_{ \rm loc}(\R^{3})$, and therefore $ \bu \in \bC^ \infty(\R^{3})$ and $ p\in C^{\infty}(\R^{3})$. 
In particular, we have for all $ \zeta \in C^{\infty}_{ \rm c}(\R^{3})$ 
\begin{align}
\intl_{ \R^{3}} | \nabla \bu |^2 \zeta dx &= \frac{1}{2}  \intl_{ \R^{3}} |\bu |^2 \Delta  \zeta dx
 + \frac{1}{2}\intl_{ \R^{3}} | \bu |^2 \bu \cdot \nabla \zeta dx+ \intl_{ \R^{3}} p \bu \cdot \nabla \zeta dx. 
\label{4.1}
\end{align}

On the basis of \eqref{4.1} we have the following Caccioppoli-type inequality. 
\begin{equation}
\intl_{B_R} | \nabla \bu | ^2 dx  
 \lesssim  R^{ -1} \bigg\{1+ \intl_{B_{ 2R}} | \bu |^3 dx\bigg\}.
\label{4.1a}
\end{equation}
  
\noindent{\bf Proof of \eqref{4.1a}}: Let $ R \le  r < \rho \le 2R$. Into \eqref{4.1} we insert a off function $ \zeta \in C^{\infty}_{\rm c}(B_\rho )$ such that $ \zeta \equiv 1$ on $ B_r$, $ 0 \le \zeta \le 1$ in $ \R^{3}$ and $ | \nabla \zeta| ^2+ | \nabla ^2\zeta |
  \lesssim  (\rho -r)^{ -2}$.  This together with H\"older's inequality and Young's inequality immediately gives \begin{align}
 \intl_{B_r} | \nabla \bu | ^2  dx &\lesssim  (\rho -r)^{ -2} \intl_{B_{\rho }} | \bu |^2dx  +   (\rho -r)^{ -1} \intl_{B_{ \rho }} | \bu |^3 dx
 + (\rho -r)^{ -1} \intl_{B_\rho }   | p- p_{ B_{ \rho }}| \, | \bu |  dx
  \cr
   &\lesssim  (\rho -r)^{ -1} \bigg\{1+\intl_{B_{ 2R} } | \bu |^3dx \bigg\}+ (\rho -r)^{ -1} \intl_{B_\rho }   | p- p_{ B_{ \rho }}| \, | \bu | dx. 
  \label{4.1b}
  \end{align}
Using H\"older's inequality,, Young's inequality and consulting  Theorem III.3.1, Theorem III.5.2 of \cite{Galdi}, we estimate the last integral involving the pressure as follows
\begin{align*}
& (\rho -r)^{ -1}\intl_{B_\rho } 
 | p- p_{ B_\rho }|\, | \bu|   dx
 \\
 &\qquad \qquad \lesssim  (\rho -r)^{ -1}  
 \bigg(\intl_{B_\rho }  | \nabla \bu|^{ \frac23} dx+ \intl_{B_\rho }| \bu |^3 dx\bigg)^{ \frac23}
  \bigg(\intl_{B_\rho } | \bu |^3 dx \bigg)^{ \frac13}
\\
 &\qquad \qquad \lesssim  \rho ^{ 1/2}(\rho -r)^{ -1}  \bigg(\intl_{B_\rho } | \nabla \bu |^2 dx \bigg)^{ \frac12}
 \bigg(\intl_{B_{ 2R}} | \bu |^3dx\bigg)^{ \frac13}
 + (\rho -r)^{ -1} \intl_{B_{ 2R} } | \bu |^3 dx
\\
 &\qquad \qquad \lesssim  \delta \intl_{B_\rho } | \nabla \bu |^2 dx + \rho (\rho -r)^{ -2} \bigg( \intl_{B_{ 2R} } | \bu |^3dx\bigg)^{ \frac23}
 + (\rho -r)^{ -1} \intl_{B_{ 2R} } | \bu |^3 dx
\\
 &\qquad \qquad \lesssim  \delta \intl_{B_\rho } | \nabla \bu |^2 dx
 + (\rho -r)^{ -1} \bigg\{1+\intl_{B_{ 2R} } | \bu |^3dx\bigg\}.
  \end{align*}
Inserting this inequality into the right-hand side of \eqref{4.1b}, we arrive at 
 \begin{equation}
\intl_{B_r} | \nabla \bu | ^2  dx  \lesssim  (\rho -r)^{ -1} \bigg\{1+ \intl_{B_{ 2R}} | \bu |^3 dx\bigg\}
+ \delta \intl_{B_\rho } | \nabla \bu | ^2 dx.
\label{4.1c}
\end{equation}
In \eqref{4.1c} taking $ \delta >0$ sufficiently small, and applying a well known iteration argument, we obtain \eqref{4.1a}.  
This completes the proof of  \eqref{4.1a}.  \hfill \Beweisende \\

\noindent{\bf Proof of Theorem 1.1 } 
Let $ \var >0$ be arbitrarily chosen, but fixed. Thanks to Lemma\,\ref{lem2.3}, in view of \eqref{3.9a} we may choose  $ R \ge  \var ^{ -1}$ such that 
\begin{equation}
\intl_{B_{ R}  \setminus B_{ R/2}}  \phi_{ \frac92} (| \bu |) dx \le  \frac{\var }{\log R}. 
\label{4.2}
\end{equation} 
Let $ \zeta \in C^{\infty}_{\rm c} (B_R)$ be a cut off function such that $ 0 \le \zeta  \le 1$ in $ B_R$, 
$ \zeta \equiv 1$ on $ B_{ R/2}$, and $ | \nabla \zeta |  \lesssim  R^{ -1}, | \nabla ^2 \zeta |  \lesssim  R^{ -2} $. 
Then from \eqref{4.1} we deduce  
\begin{equation}
\intl_{B_{ R/2}} | \nabla \bu |^2 dx \lesssim  R^{ -2}  \intl_{B_R  \setminus B_{ R/2}} |\bu |^2 dx
 + R^{ -1}\intl_{B_R  \setminus B_{ R/2}} | \bu |^3 dx + R^{ -1}\intl_{B_R  \setminus B_{ R/2}} 
 | p- p_{ B_R  \setminus B_{ R/2}}|\, | \bu| dx
\label{4.3}
\end{equation}
Using H\"older's inequality, Young's inequality and consulting \cite{Galdi}, we estimate the last integral involving the pressure as follows
\begin{align*}
& R^{ -1}\intl_{B_R  \setminus B_{ R/2}} 
 | p- p_{ B_R  \setminus B_{ R/2}}|\, | \bu|   dx
 \\
 &\qquad \qquad \lesssim  R^{ -1}  
 \bigg(\intl_{B_R  \setminus B_{ R/2}}  | \nabla \bu|^{ \frac32}dx  +\intl_{B_R  \setminus B_{ R/2}} | \bu |^3 dx\bigg)^{\frac23}
  \bigg(\intl_{B_R  \setminus B_{ R/2}} | \bu |^3dx\bigg)^{ \frac13}
\\
 &\qquad \qquad \lesssim  R^{ -\frac12} \bigg(\intl_{B_R} | \nabla \bu |^2dx \bigg)^{ \frac12} \bigg(\intl_{B_R  \setminus B_{ R/2}} | \bu |^3dx \bigg)^{ \frac13}
 + R^{ -1} \intl_{B_R  \setminus B_{ R/2}} | \bu |^3 dx
\\
 &\qquad \qquad \lesssim  R^{ -\frac13} \bigg(\intl_{B_R} | \nabla \bu |^2dx \bigg)^{ \frac34} 
+ R^{ -1} \intl_{B_R  \setminus B_{ R/2}} | \bu |^3dx
\\
 &\qquad \qquad \lesssim  R^{ -\frac18}+  R^{ -\frac16}\intl_{B_R} | \nabla \bu |^2dx + 
R^{ -1} \intl_{B_R  \setminus B_{ R/2}} | \bu |^3dx . 
 \end{align*}
 Once more using H\"older's inequality along with Young's inequality we easily find 
\[
R^{ -2}  \intl_{B_R  \setminus B_{ R/2}} |\bu |^2  dx \le R^{ -1}+ R^{ -1}  \intl_{B_R  \setminus B_{ R/2}} |\bu |^3 dx.  
\] 
Inserting the last two  inequalities into the right-hand side of \eqref{4.3}, we arrive at
\begin{equation}
\intl_{B_{ R/2}} | \nabla \bu |^2 dx  \lesssim   R^{ -1}\intl_{B_R  \setminus B_{ R/2}} | \bu |^3dx  + R^{ -\frac18}
+  R^{ -\frac16} \intl_{B_R} | \nabla \bu |^2 dx. 
\label{4.4}
\end{equation}
We now estimate the last integral on the right-hand side of \eqref{4.4} by means of \eqref{4.1a}. This implies 
\begin{equation}
\intl_{B_{ R/2}} | \nabla \bu |^2  dx \lesssim   R^{ -1}\intl_{B_R  \setminus B_{ R/2}} | \bu |^3 dx + R^{ -\frac18}
+  R^{ -\frac76} \intl_{B_{ 2R}} | \bu |^3 dx. 
\label{4.5}
\end{equation}  
By our assumption \eqref{2.1aa} we know that $ \bu \in L^q(\R^{3})$ for all $ q > \frac{9}{2}$. 
This follows from standard regularity theory of the steady Navier-Stokes equations (e.g. see \cite{ser2}). 
For $ \frac{9}{2} <q < \frac{54}{11}$ 
we find with  the help of Jensen's inequality  
\begin{equation}
R^{ -\frac18}+R^{ -\frac76} \intl_{B_{ 2R}} | \bu |^3 dx \lesssim R^{ -\frac18}+  R ^{ \frac{3q-9}{q}- \frac{7}{6}} \| \bu \|_{ q} \rightarrow 0\quad  
\text{ as}\quad  R \rightarrow +\infty. 
\label{4.9}
\end{equation}

Noting that $ \phi _{ 3/2}$ is convex, applying Jensen's inequality,  we get 
\begin{align*}
\lefteqn{\phi_{\frac32} \bigg(\frac{8}{7R\mes ( B_1)} \intl_{B_{ R}  \setminus B_{ R/2}}  | \bu |^3 dx\bigg) =\phi_{\frac32} \bigg( R^2  \intmw_{B_{ R}  \setminus B_{ R/2}}   |v|^3 dx \bigg)}
\hspace{1.in}\\
&\leq
 \intmw_{B_{ R}  \setminus B_{ R/2}}  \phi _{ \frac32}(R^2| \bu |^3) dx
\\
 &\lesssim   \int_{B_{ R}  \setminus B_{ R/2}} \frac{| \bu |^{ \frac92}}{\log \frac{1+ 2R^2 | \bu |^3}{R^2 | \bu |^3}} dx.
 \end{align*}
We split the integral on the right-hand side into two parts by setting 
\begin{align*}
A_1 &= \{x\in B_{ R}  \setminus B_{ R/2} \,|\, | \bu |^3 \le \var R^{ -2}\},\quad  
\\
A_2 &= \{x\in B_{ R}  \setminus B_{ R/2} \,|\, | \bu |^3 >\var R^{ -2}\}.
\end{align*}

Firstly, we easily see that 
\[
\int_{A_1} \frac{| \bu |^{ \frac92}}{\log \frac{1+ 2R^2 | \bu |^3}{R^2 | \bu |^3}}  dx \lesssim  \var ^{ \frac32}.
\]
Secondly, with help of Lemma\,\ref{lem3.4} and recalling that $ R \ge \frac{1}{\var }$ we have in  $ A_2$
\begin{align*}
4 \log R &\ge  \log R^2 + \log \frac{1}{\var } + \log 2 
= \log 2 \frac{R^2}{\var } \ge \log \frac{1 + 2 \var R^{ -2}}{\var R^{ -2}} 
\\
&\ge \log\frac{1+ 2| \bu |^3}{| \bu |^3} \gtrsim \log\frac{1+ 2| \bu |}{| \bu |}.
\end{align*}
With this estimate along with \eqref{4.2} we get
\begin{align*}
\int_{A_2} \frac{| \bu |^{ \frac92}}{\log \frac{1+ 2R^2 | \bu |^3}{R^2 | \bu |^3}} dx &\lesssim  \frac{1}{\log 2}\intl_{A_2} | \bu |^{ \frac92} dx \lesssim    \frac{\log R}{\log 2} \int_{B_{ R}  \setminus B_{ R/2}} \frac{| \bu |^{ \frac92}}{\log \frac{1+2| \bu |}{| \bu |}} dx
\\
& \lesssim   \log R \int_{B_{ R}  \setminus B_{ R/2}} \phi _{ \frac92} (| \bu |)  dx \lesssim \var. 
\end{align*}
Accordingly,
\[
 \phi_{ \frac32} \bigg(\frac{8}{7R \mes (B_1)} \intl_{B_{ R}  \setminus B_{ R/2}}  | \bu |^3 dx\bigg)  \lesssim  \var.  
\]
Thus, in view \eqref{4.5} together with the estimates we have just obtained we are able to chose a sequence $ R_k \rightarrow +\infty$ as $ k \rightarrow +\infty$, such that 
\[
\intl_{B_{ R_k/2}} | \nabla \bu |^2 dx \rightarrow 0\quad  \text{ {\it as}}\quad  k \rightarrow +\infty,
\]
which yields $ \nabla \bu ={ 0}$ and therefore $ \bu \equiv \const={0}$.   \hfill \Beweisende

\subsection {Proof of Theorem 1.2}
   \noindent{\bf Proof  of  Theorem 1.2 } Let us denote the vorticity $\o=\mathrm{curl}\, v$. Then, it is well-known that from (NS)  that the following equation holds true.
\bb\label{21}
\Delta Q-v\cdot\nabla Q=|\o|^2.
\ee
Under the condition (\ref{13}) we have $Q(x)\leq 0$ for all $x\in \Bbb R^3$ by the maximum principle applied to (\ref{21}). 
Moreover, by the maximum principle again, either $Q(x)\equiv 0$ on $\Bbb R^3$, or $Q(x)<0$ for all $x\in \Bbb R^3$.
Indeed, any point $x_0\in \Bbb R^3$ such that $Q(x_0)=0$ is a point of local maximum, which is not allowed unless $Q\equiv 0$ by the maximum principle.
 Let  $Q(x)\not\equiv0$ on $\Bbb R^3$,  then without the loss of generality we may assume $|Q(x)| >0$ for all $x\in \Bbb R^3$.
We  set $\sup_{x\in \Bbb R^3} |Q|=M>0$.  
Let $f\in C(\Bbb R)$. For $\l \in [0, M)$ we set $D_\l =\{ x\in \Bbb R^3\, |\, |Q(x)| >\l \}$. Then, we compute
\bq\label{23}
\lefteqn{\int_{D_\l }  f(Q(x))v\cdot \nabla Q \, dx=
\int_{D_\l } v\cdot \nabla \left(\int_0 ^{Q(x)} f(q) dq\right) dx}\n \\
&&= \int_{  D_\l} \mathrm{div}\, \left(v \int_0 ^{Q(x)} f(q ) dq \right) dx = \int_{\partial D_\l}  \left(\int_0 ^{Q(x)} f(q) dq\right)\, v\cdot \nu\, dS\n \\
&&= \int_0 ^{ \l }  f(q) dq \int_{\partial D_\l}   v\cdot \nu \, dS= \int_0 ^{ \l }  f(q) dq \int_{D_\l}   \mathrm{div} \, v\, dx=0.
\eq 
where $\nu=\nabla Q/|\nabla Q|$ is the outward unit normal vector on  $\partial D_\l$.
For $\l \in (0, M)$ we Integrate (\ref{21}) over $D_\l$. Then, using the fact (\ref{23}), we have
\bq\label{24}
\int_{D_\l} |\o|^2 dx &=& \int_{D_\l} \Delta Q \,dx =\int_{\partial D_\l} \frac{\partial Q}{\partial \nu}  dS\n \\
&=& \int_{\partial D_\l} |\nabla Q| dS.
\eq
Using the co-area formula, we obtain
\bqn
\lefteqn{\int_{D_\l}\frac{|\nabla Q|^2}{|Q|} \left( \log \frac{eM }{|Q|}\right)^{-\a-1}  dx=\int_\l ^M \int_{\partial D_q } \frac{|\nabla Q|}{|Q|} \left( \log \frac{eM }{|Q|}\right)^{-\a-1} dS dq}\hspace{1.in}\n \\
&&=\int_\l ^M\frac{1}{q} \left( \log \frac{eM }{q}\right)^{-\a-1} \int_{\partial D_q} |\nabla Q| dS dq\n \\
&&\leq \int_\l ^M\frac{1}{q} \left( \log \frac{eM }{q}\right)^{-\a-1} dq \int_{\partial D_\l} |\nabla Q| dS \n \\
&&=\frac{1}{\a}\left\{1- \left( \log \frac{eM }{\l}\right)^{-\a}\right\} \int_{D_\l} |\o|^2 dx\n \\
&&\leq \frac{1}{\a} \int_{\Bbb R^3} |\o|^2 dx,\eqn
where we used (\ref{24}) in the fourth line.
Passing $\l \to 0$, and applying the monotone convergence theorem, we obtain (\ref{th11}).
Next, we assume (\ref{th12}) holds.
We consider  a standard cut-off function $\sigma\in C_0
^\infty([0, \infty))$ such that $ \sigma(s)=1$ if $s<1$, and $\s(s)=0$ if $s>2$, 
and $0\leq \sigma  (s)\leq 1$ for $1<s<2$. 
For each $\a \in (0, 1)$  we  define $\s_\a(x):=\s_\a (Q(x))\in C_0 ^\infty (\Bbb R^3)$ by
$$
\s_\a (x)= 1-\s \left\{  3\left(\log \frac{eM}{|Q(x)|} \right)^{-\a} \right\}.
$$
We note that
$$
\left\{ \aligned & \s_\a (x)=1,  \quad \mbox{if}\quad   |Q(x)|\geq  Me^{1- \left(\frac32\right)^{\frac{1}{\a}} }  ,\\
 &0<\s_\a (x)<1,  \quad \mbox{if}\quad Me^{1-3 ^{\frac{1}{\a}}}< |Q(x)| <Me^{1- \left(\frac32\right)^{\frac{1}{\a}} } ,\\
 & \s _\a(x)=0,\quad \mbox{if}\quad   |Q(x)| \leq Me^{1-3 ^{\frac{1}{\a}}}.\endaligned \right.
 $$
We multiply (\ref{21}) by $\s_\a$, and integrate it over $\Bbb R^3$.  Then, the convection term vanishes by (\ref{23}).
Let $\a_1 >0$ be fixed. For all $\a >\a_1$ we have 
\bqn
\lefteqn{ \int_{\Bbb R^3} |\o|^2 \s_{\a_1} (x)dx\leq \int_{\Bbb R^3} |\o|^2 \s_{\a} (x)dx=
 \int_{\Bbb R^3} \Delta Q\, \s_\a (x) dx}\n \\
&=& -3\a \int_{\{ Me^{1-3 ^{\frac{1}{\a}}}< |Q(x)| <Me^{1- \left(\frac32\right)^{\frac{1}{\a}} }\}  } \frac{|\nabla Q|^2}{|Q|} \left( \log \frac{eM}{|Q|} 
\right)^{-\a -1} \s' \left\{ 3\left( \log \frac{eM}{|Q(x)|} \right)^{-\a}\right\} dx\n \\
&\leq& 3\a \sup_{1\leq s\leq 2} |\s'(s)| \int_{\{  Me^{1-3 ^{\frac{1}{\a}}}< |Q(x)| <Me^{1- \left(\frac32\right)^{\frac{1}{\a}} }  \}} \frac{|\nabla Q|^2}{|Q|} \left( \log \frac{eM}{|Q|} \right)^{-\a -1}  dx\n \\
&\to &0 \quad \mbox{as}\quad \a\to 0.
\eqn
Hence, we have shown
$ \int_{\Bbb R^3} |\o|^2 \s_{\a_1}(x)dx= 0
$
for all $\a_1 >0$, which implies that  $\o=0$ on $\Bbb R^3$. This, combined with the fact div $v=0$ implies that $v$ is a harmonic function on $\Bbb R^3$.
The boundary condition, together with the Liouville theorem for harmonic function, leads us to conclude $v=0$ on $\Bbb R^3$.  \hfill \Beweisende  \\

\subsection{Proof of Theorem 1.3}

We first establish integrability conditions on the vector fields for the Biot-Savart's formula in $\mathbb{R}^3$.
\begin{pro}\label{pro1}
Let $\xi=\xi(x)=(\xi_1(x),\xi_2(x),\xi_3(x))$ and $\eta=\eta(x)=(\eta_1(x),\eta_2(x),\eta_3(x))$ be smooth vector fields on $\mathbb{R}^3$.
Suppose there exists $q\in [1,3)$ such that $\eta \in L^q(\mathbb{R}^3)$. Let $\xi$ solve
\begin{equation}
\Delta\xi=-\nabla\times\eta, \label{2.1}
\end{equation}
under the boundary condition; either
\begin{equation}
|\xi(x)|\rightarrow0\quad\text{ as}\quad|x|\rightarrow+\infty, \label{2.2}
\end{equation}
or
\begin{equation}
\xi \in L^{s}(\mathbb{R}^3)\quad\text{for some}\quad s\in [1,\infty).\label{2.3}
\end{equation}
Then, the solution of \eqref{2.1} is given by
\begin{equation}
\xi(x)=\frac{1}{4\pi}\int_{\mathbb{R}^3}\frac{(x-y)\times\eta(y)}{|x-y|^3} dy\quad \forall x\in \mathbb{R}^3. \label{2.4}
\end{equation}
\end{pro}

\noindent{\bf Proof} \, Let  $\sigma \in C^{\infty}_{0}(\mathbb{R}^3)$  be the cut-off function defined in the proof of Theorem 1.1. For each $R>0$ we define $\sigma_{R}(x):=\sigma\left(\frac{|x|}{R}\right)$. Given $\epsilon>0$
we denote $B_{\epsilon}(y)=\{x\in\mathbb{R}^3||x-y|<\epsilon\}$. Let us fix $y\in\mathbb{R}^3$ and $\epsilon\in(0,\frac{R}{2})$.
We multiply \eqref{2.1} by $\frac{\sigma_{R}(|x-y)|}{|x-y|}$, and integrate it with respect to the variable $x$ over $\mathbb{R}^3 \setminus B_{\epsilon}(y)$.
Then,
\begin{equation}
\int_{\{|x-y|>\epsilon\}}\frac{\Delta\xi \sigma_{R}}{|x-y|} dx= -\int_{\{|x-y|>\epsilon\}} \frac{\sigma_{R}\nabla\times\eta(y)}{|x-y|}dx. \label{2.5}
\end{equation}
Since $\Delta \frac{1}{|x-y|}=0$ on $\mathbb{R}^3\setminus B_{\epsilon}(y)$, one has
\begin{align*}
\frac{\Delta\xi \sigma_{R}}{|x-y|}&=\sum^{3}_{i=0} \partial_{x_i}\left(\frac{\partial_{x_i}\xi\sigma_{R}}{|x-y|}\right)-
\sum^{3}_{i=0} \partial_{x_i}\left(\frac{\xi\partial_{x_i}\sigma_{R}}{|x-y|}\right)\\
&-\sum^{3}_{i=0}\partial_{x_i}\left(\xi\sigma_{R}\partial_{x_i}\left(\frac{1}{|x-y|}\right)\right)+\frac{\xi\Delta\sigma_{R}}{|x-y|}
+2\sum^{3}_{i=0}\xi\partial_{x_i}\left(\frac{1}{|x-y|}\right)\partial_{x_i}\sigma_{R}.
\end{align*}
Therefore, applying the divergence theorem, and observing $\partial_{\nu}\sigma_{R} = 0 $ on $\partial B_{\epsilon}(y)$,
we have
\begin{align}\label{2.6}
\int_{\{|x-y|>\epsilon\}} &\frac{\Delta\xi\sigma_{R}}{|x-y|} dx = \int_{\{|x-y|=\epsilon\}}\frac{\partial_{\nu}\xi}{|x-y|}dS\nonumber\\
&-\int_{\{|x-y|=\epsilon\}}\frac{\xi}{|x-y|^2}dS+\int_{\{|x-y|>\epsilon\}}\frac{\xi\Delta\sigma_{R}}{|x-y|}dS\nonumber\\
&-2\int_{\{|x-y|>\epsilon\}}\frac{(x-y)\cdot\nabla \sigma_{R} \xi}{|x-y|^3}dS
\end{align}
where $\partial_{\nu}(\cdot)$ denotes the outward normal derivative on $\partial B_{\epsilon}(y)$. Passing $\epsilon\rightarrow0$,
one can easily compute that
\begin{align}\label{2.7}
\text{RHS of \eqref{2.6}}\quad&\rightarrow\quad -4\pi\xi(y) + \int_{\mathbb{R}^3} \frac{\xi\Delta\sigma_{R}}{|x-y|}dx
-2\int_{\mathbb{R}^3}\frac{(x-y)\cdot\nabla \sigma_{R} \xi}{|x-y|^3}dx\nonumber\\
&:= \quad I_1 + I_2 +I_3
\end{align}
Next, using the formula
$$\frac{\s_{R}\nabla\times\eta}{|x-y|}=\nabla\times\left(\frac{\s_{R}\eta}{|x-y|}\right)-\frac{\nabla\s_R \times \eta}{|x-y|}
+\frac{(x-y)\times\eta\s_R}{|x-y|^3},$$
and using the divergence theorem, we obtain the following representation for the right hand side of \eqref{2.5}.
\begin{align}\label{2.8}
&\int_{\{|x-y|>\epsilon\}} \frac{\s_R\nabla\times\eta}{|x-y|}dx= \int_{\{|x-y|=\epsilon\}}\nu\times\left(\frac{\eta}{|x-y|}\right)dS\nonumber\\
&\quad -\int_{\{|x-y|>\epsilon\}}\frac{\nabla\s_R\times\eta}{|x-y|}dx+\int_{\{|x-y|>\epsilon\}} \frac{(x-y)\times\eta\s_R}{|x-y|^3}dx,
\end{align}
where we denoted $\nu=\frac{y-x}{|y-x|}$, the outward unit normal vector on $\partial B_{\epsilon}(y)$. Passing $\epsilon\rightarrow0$,
we easily deduce
\begin{align}\label{2.9}
\text{RHS of \eqref{2.8}}\quad&\rightarrow \quad- \int_{\mathbb{R}^3} \frac{\nabla\s_R\times\eta}{|x-y|}dx
-2\int_{\mathbb{R}^3}\frac{(x-y)\times\eta\sigma_{R}}{|x-y|^3}dx\nonumber\\
&:= \quad J_1 + J_2  \quad\text{as}\quad\epsilon\rightarrow0.
\end{align}
We now pass $R\rightarrow\infty$ for each term of \eqref{2.7} and \eqref{2.9} respectively below. Under the boundary
condition \eqref{2.2} we estimate:
\begin{align*}
|I_2|\,\,&\leq\,\, \int_{\{R\leq|x-y|\leq2R\}}\frac{|\xi(x)||\Delta\s_R(x-y)|}{|x-y|}dx\\
         &\leq\,\,\frac{\|\Delta\s\|_{L^{\infty}}}{R^2}\sup_{R\leq|x-y|\leq2R}|\xi(x)|
         \left(\int_{\{R\leq|x-y|\leq2R\}}dx\right)^{\frac{2}{3}}\left(\int_{\{R\leq|x-y|\leq2R\}}\frac{dx}{|x-y|^3}\right)^{\frac{1}{3}}\\
         & \lesssim \,\, \|\Delta\s\|_{L^{\infty}}\left(\int_{R}^{2R}\frac{dr}{r}\right)^{\frac{1}{3}}\sup_{R\leq|x-y|\leq2R}|\xi(x)|\rightarrow0
\end{align*}
as $R\rightarrow\infty$ by the assumption \eqref{2.2}, while under the condition \eqref{2.3} we have
\begin{align*}
|I_2|\,\,&\leq\,\, \int_{\{R\leq|x-y|\leq2R\}}\frac{|\xi(x)||\Delta\s_R(x-y)|}{|x-y|}dx\\
          &\leq\,\,\frac{\|\Delta\s\|_{L^{\infty}}}{R^2}\|\xi\|_{L^s} \left(\int_{\{R\leq|x-y|\leq2R\}}\frac{dx}{|x-y|^{\frac{s}{s-1}}}\right)^{\frac{s-1}{s}}\\
         & \lesssim \,\, R^{-\frac{3}{s}} \|\Delta\s\|_{L^{\infty}}\|\xi\|_{L^s}\rightarrow0
\end{align*}
as $R\rightarrow\infty$. Similarly, under \eqref{2.2}
\begin{align*}
|I_3|\,\,&\leq\,\, 2\int_{\{R\leq|x-y|\leq2R\}}\frac{|\xi(x)||\nabla\s_R(x-y)|}{|x-y|^2}dx\\
         & \lesssim \,\,\frac{\|\nabla\s\|_{L^{\infty}}}{R}\sup_{R\leq|x-y|\leq2R}|\xi(x)|
         \left(\int_{\{R\leq|x-y|\leq2R\}}dx\right)^{\frac{1}{3}}\left(\int_{\{R\leq|x-y|\leq2R\}}\frac{dx}{|x-y|^3}\right)^{\frac{2}{3}}\\
         & \lesssim \,\, \|\nabla\s\|_{L^{\infty}}\left(\int_{R}^{2R}\frac{dr}{r}\right)^{\frac{2}{3}}\sup_{R\leq|x-y|\leq2R}|\xi(x)|\rightarrow0
\end{align*}
as $R\rightarrow\infty$, while under the condition \eqref{2.3} we estimate
\begin{align*}
|I_3|\,\,&\leq\,\, 2\int_{\{R\leq|x-y|\leq2R\}}\frac{|\xi(x)||\nabla\s_R(x-y)|}{|x-y|^2}dx\\
          & \lesssim \,\,\frac{\|\nabla\s\|_{L^{\infty}}}{R^2}\|\xi\|_{L^s} \left(\int_{\{R\leq|x-y|\leq2R\}}\frac{dx}{|x-y|^{\frac{2s}{s-1}}}\right)^{\frac{s-1}{s}}\\
         & \lesssim \,\, R^{-\frac{3}{s}} \|\nabla\s\|_{L^{\infty}}\|\xi\|_{L^s}\rightarrow0
\end{align*}
as $R\rightarrow\infty$. Therefore, the right hand side of \eqref{2.6} converges to $-4\pi\xi(y)$ as as $R\rightarrow\infty$.
For $J_1,J_2$ we estimate
\begin{align*}
|J_1|\,\,&\leq\,\, \int_{\{R\leq|x-y|\leq2R\}}\frac{|\nabla\s_R||\eta|}{|x-y|}dx\\
          &\leq\,\,\frac{\|\nabla\s\|_{L^{\infty}}}{R}\|\eta\|_{{L^q}(R\leq|x-y|\leq2R)} \left(\int_{\{R\leq|x-y|\leq2R\}}\frac{dx}{|x-y|^{\frac{q}{q-1}}}\right)^{\frac{q-1}{q}}\\
         & \lesssim \,\, \|\nabla\s\|_{L^{\infty}}\|\eta\|_{{L^q}(R\leq|x-y|\leq2R)}R^{-\frac{2}{q}}\rightarrow0
\end{align*}
as $R\rightarrow\infty$. In passing $R\rightarrow\infty$ in $J_2$ of \eqref{2.9}, in order to use the dominated convergence theorem, we estimate
\begin{align}\label{2.10}
\int_{\mathbb{R}^3}\left|\frac{(x-y)\times\eta(y)}{|x-y|^3}\right|dx\,\,
&\leq\,\,\int_{\{|x-y|<1\}}\frac{|\eta|}{|x-y|^2}dx+\int_{\{|x-y|\geq1\}}\frac{|\eta|}{|x-y|^2}dx\nonumber\\
&:=\,\, J_{21}+J_{22}.
\end{align}
$J_{21}$ is easy to handle as follows.
\begin{equation}\label{2.11}
J_{21}\leq\|\eta\|_{L^{\infty}(B_1(y))}\int_{\{|x-y|<1\}}\frac{dx}{|x-y|^2}=4\pi\|\eta\|_{L^{\infty}(B_1(y))}<+\infty
\end{equation}
For $J_{22}$ we estimate
\begin{align}\label{2.12}
J_{22}\,\,&\leq\,\,\left(\int_{\mathbb{R}^3}|\eta|^q\right)^{\frac{1}{q}}
  \left(\int_{\{|x-y|>1\}}\frac{dx}{|x-y|^{\frac{2q}{q-1}}}\right)^{\frac{q-1}{q}}\nonumber\\
  & \lesssim \,\,\|\eta\|_{L^q}\left(\int_{1}^{\infty}r^{\frac{-2}{q-1}}dr\right)^{\frac{q-1}{q}}<+\infty,
\end{align}
if $1<q<3$. In the case of $q=1$ we estimate simply
\begin{equation}\label{2.13}
J_{22}\leq\int_{\{|x-y|>1\}}|\eta|dx \leq \|\eta\|_{L^1}.
\end{equation}
Estimates of \eqref{2.10}-\eqref{2.13} imply
$$\int_{\mathbb{R}^3}\left|\frac{(x-y)\times\eta(y)}{|x-y|^3}\right| dx<+\infty.$$
Summarising the above computations, one can pass first $\epsilon\rightarrow0$, and then $R\rightarrow+\infty$ in \eqref{2.5},
applying the dominated convergence theorem, to obtain finally \eqref{2.4}.  \hfill \Beweisende \\

\begin{cor}\label{cor1}
Let $v$ be a smooth solution to (NS) satisfying (\ref{11})  such that $\w \in L^q(\mathbb{R}^3)$ for some $q\in[\frac{3}{2},3)$.
Then, we have
\begin{equation}\label{2.14}
v(x)=\frac{1}{4\pi}\int_{\mathbb{R}^3}\frac{(x-y)\times \w(y)}{|x-y|^3}dy,
\end{equation}
and
\begin{equation}\label{2.15}
\w(x)=\frac{1}{4\pi}\int_{\mathbb{R}^3}\frac{(x-y)\times (v(y)\times \w(y))}{|x-y|^3}dy.
\end{equation}
\end{cor}
\noindent{\bf Proof} \, Taking curl of the defining equation of the vorticity, $\nabla\times v=\w$, using $\text{div}\,v=0$, we have
$$\Delta v=-\nabla \times \w ,$$
which provides us with \eqref{2.14} immediately by application of Proposition \ref{pro1}.
In order to show \eqref{2.15} we recall that, using the vector identity $\frac{1}{2}\nabla|v|^2=(v\cdot\nabla)v+v\times(\nabla\times v)$,
one can rewrite (NS)  as
$$-v\times \w= -\nabla\left(p+\frac{1}{2}|v|^2\right) +\Delta v.$$
Taking curl on this, we obtain
$$\Delta \w = -\nabla \times(v\times \w).$$
The formula \eqref{2.15} is deduced immediately from this equations by applying the proposition \ref{pro1}.
For the allowed rage of $q$ we recall the Sobolev and the Calderon-Zygmund inequalities(\cite{Stein}),
\begin{equation}\label{2.16}
\|v\|_{L^{\frac{3q}{3-q}}}  \lesssim  \|\nabla v\|_{L^q} \lesssim \|\w\|_{L^q},\quad 1<q<3,
\end{equation}
which imply $v \times \w \in L^{\frac{3q}{6-q}}(\mathbb{R}^3)$ if $ \w\in L^{q}(\mathbb{R}^3)$. We also note that $\frac{3}{2}\leq q <3 $
if and only if $1\leq\frac{3q}{6-q}<3$.  \hfill \Beweisende \\

\noindent{\bf Proof of Theorem 1.3  } Under the hypothesis \eqref{11} and $\w\in L^q(\mathbb{R}^3)$ with $q \in [\frac{3}{2},3)$
both of the relations \eqref{2.14} and \eqref{2.15} are valid. We first prove the following.\\

\noindent{\it{Claim: }} For each $x,y \in \mathbb{R}^3$
\begin{equation}\label{2.17}
0\leq |\w(x)|^2 = \int_{\mathbb{R}^3} \Phi(x,y) dy \leq \int_{\mathbb{R}^3} |\Phi(x,y)|dy < +\infty,
\end{equation}
and
\begin{equation}\label{2.18}
0= \int_{\mathbb{R}^3} \Phi(x,y) dx \leq \int_{\mathbb{R}^3} |\Phi(x,y)|dx < +\infty.
\end{equation}
\noindent{\bf Proof of }(\ref{2.19}):  Decomposing the integral and using the H\"{o}lder inequality, we estimate
\begin{align}\label{2.20}
\int_{\mathbb{R}^3}|\Phi(x,y)|dy\,\,&\leq\,\,|\w(x)|\left(\int_{\{|x-y|\leq1 \}}\frac{|v(y)||\w(y)|}{|x-y|^2}dy
    +\int_{\{|x-y|>1 \}}\frac{|v(y)||\w(y)|}{|x-y|^2}dy \right)\nonumber \\
    &\leq\,\,|\w(x)|\|v\|_{L^{\infty}(B_1(x))}\|\w\|_{L^{\infty}(B_1(x))}\int_{\{|x-y|\leq1 \}}\frac{dy}{|x-y|^2}\nonumber\\
     &\qquad \quad +|\w(x)|\|v\|_{L^{\frac{3q}{3-q}}}\|\w\|_{L^q}\left(\int_{\{|x-y|\geq1 \}}\frac{dy}{|x-y|^{\frac{6q}{4q-6}}}\right)^{\frac{4q-6}{3q}}\nonumber\\
 & \lesssim \,\,|\w(x)|\|v\|_{L^{\infty}(B_1(x))}\|\w\|_{L^{\infty}(B_1(x))}\nonumber\\
     &\qquad\quad +|\w(x)|\|\w\|^2_{L^q}\left(\int_{1}^{\infty}r^{\frac{q-6}{2q-3}}dr\right)^{\frac{4q-6}{3q}} <+\infty,
    \end{align}
 where we used \eqref{2.16} and the fact that $\frac{q-6}{3q-3}<-1$ if $\frac{3}{2}<q<3$. In the case $q=\frac{3}{2}$ we estimate, instead,
\begin{align}\label{2.21}
\int_{\mathbb{R}^3}|\Phi(x,y)|dy\,\,&\leq\,\,|\w(x)|\left(\int_{\{|x-y|\leq1 \}}\frac{|v(y)||\w(y)|}{|x-y|^2}dy
    +\int_{\{|x-y|>1 \}}\frac{|v(y)||\w(y)|}{|x-y|^2}dy \right)\nonumber \\
    &\leq\,\,|\w(x)|\|v\|_{L^{\infty}(B_1(x))}\|\w\|_{L^{\infty}(B_1(x))}+|\w(x)|\|v\|_{L^3}\|\w\|_{L^{\frac{3}{2}}}<+\infty.
\end{align}
We also have
\begin{align}\label{2.22}
\int_{\mathbb{R}^3}|&\Phi(x,y)|dx \leq |v(y)||\w(y)|\left(\int_{\{|x-y|\leq1\}}\frac{|\w(x)}{|x-y|^2}dx
       +\int_{\{|x-y|>1\}}\frac{|\w(x)}{|x-y|^2}dx\right)\nonumber\\
 & \lesssim |v(y)||\w(y)|\|\w\|_{L^{\infty}(B_1(y))} + |v(y)||\w(y)|\|\w\|_{L^q}\left(\int_{\{|x-y|>1\}}\frac{dx}{|x-y|^{\frac{2q}{q-1}}}\right)^{\frac{q-1}{q}}\nonumber\\
 & \lesssim |v(y)||\w(y)|\|\w\|_{L^{\infty}(B_1(y))}
  + |v(y)||\w(y)|\|\w\|_{L^q}\left(r^{-\frac{2}{q-1}}dr\right)^{\frac{q-1}{q}}<+\infty
\end{align}
where we used the fact that $-\frac{2}{q-1}<-1$ if $\frac{3}{2}\leq q<3$. From \eqref{2.15} we immediately obtain
\begin{align}\label{2.23}
\int_{\mathbb{R}^3}\Phi(x,y)dy\,\,&=\,\,\w(x)\cdot\left(\frac{1}{4\pi}\int_{\mathbb{R}^3}\frac{(x-y)\times (v(y)\times \w(y))}{|x-y|^3}dy\right)\nonumber\\
&=\,\,|\w(x)|^2\,\,\geq\,0,\quad \forall x\in \mathbb{R}^3
\end{align}
and combining this with \eqref{2.20}, we deduce \eqref{2.17}. On the other hand, using \eqref{2.14}, we find
\begin{align}\label{2.24}
\int_{\mathbb{R}^3}\Phi(x,y)dx\,\,&=\,\,\frac{1}{4\pi}\int_{\mathbb{R}^3}\frac{\w(x)\cdot(x-y)\times (v(y)\times \w(y))}{|x-y|^3}dx\nonumber\\
&=\,\,\left(\frac{1}{4\pi}\int_{\mathbb{R}^3}\frac{\w(x)\times (x-y)}{|x-y|^3}dx\right)\cdot v(y) \times \w(y)\nonumber\\
&=\,\, v(y)\cdot v(y) \times \w(y) = 0
\end{align}
for all $y\in \mathbb{R}^3$, and combining this with \eqref{2.22}, we have proved \eqref{2.18}. This completes the proof of the claim.\\

If (\ref{th13}) holds, then  from \eqref{2.23} and \eqref{2.24}  provide us with
$$\int_{\mathbb{R}^3}|\w(x)|^2dx=\int_{\mathbb{R}^3}\int_{\mathbb{R}^3}\Phi(x,y)dydx=\int_{\mathbb{R}^3}\int_{\mathbb{R}^3}\Phi(x,y)dxdy=0.$$
Hence,
\begin{equation}\label{2.28}
\w=0\quad\text{on}\quad{\mathbb{R}^3}.
\end{equation}
We remark parenthetically that in deriving \eqref{2.28} it is not necessary to assume that $\int_{\mathbb{R}^3}|\w(x)|^2 dx <+\infty$,
and therefore we do not need to restrict ourselves to $\w\in L^2(\mathbb{R}^3)$. Hence, from \eqref{2.14} and \eqref{2.28}, we conclude
$v=0$ on $\mathbb{R}^3$.  \hfill \Beweisende \\

\hspace{0.5cm}
$$\mbox{\bf Acknowledgements}$$
Chae was partially supported by NRF grants 2006-0093854 and  2009-0083521, while Wolf has been supported 
supported by the German Research Foundation (DFG) through the project WO1988/1-1; 612414.

\end{document}